\documentclass{amsart}
\usepackage{amsmath,amsfonts,epsfig,verbatim}

\begin{document}
\bibliographystyle{plain}

\newtheorem{thm}{Theorem}[section]
\newtheorem{lem}[thm]{Lemma}
\newtheorem{prop}[thm]{Proposition}
\newtheorem{cor}[thm]{Corollary}
\newtheorem{conj}[thm]{Conjecture}
\newtheorem{mainlem}[thm]{Main Lemma}
\newtheorem{defn}[thm]{Definition}
\newtheorem{rmk}[thm]{Remark}

\def\square{\hfill${\vcenter{\vbox{\hrule height.4pt \hbox{\vrule width.4pt
height7pt \kern7pt \vrule width.4pt} \hrule height.4pt}}}$}

\newtheorem*{namedtheorem}{\theoremname}
\newcommand{\theoremname}{testing}
\newenvironment{named}[1]{\renewcommand{\theoremname}{#1}\begin{namedtheorem}}{\end{namedtheorem}}

\newenvironment{pf}{{\it Proof:}\quad}{\square \vskip 12pt}
\newcommand{\dt}{\ensuremath{\text{det}}}
\newcommand{\U}{\ensuremath{\widetilde}}
\newcommand{\Hn}{\ensuremath{\mathbb{H}^3}}
\newcommand{\h}{{\text{hyp}}}
\newcommand{\acyl}{{\text{acyl}}}
\newcommand{\inc}{{\text{inc}}}
\newcommand{\tp}{{\text{top}}}
\newcommand{\V}{\ensuremath{\text{Vol}}}
\newcommand{\Hess}{\ensuremath{{\text{Hess} \ }}}
\title[Infinitesimal rigidity]{Infinitesimal rigidity of a compact hyperbolic $4$-orbifold with totally geodesic boundary}
\author{Tarik Aougab and Peter A. Storm}
\date{September, 2008}
\begin{abstract}
%Let $P \subset\ \mathbb{H}^{4}$ be the hyperbolic regular right-angled $120$-cell.  Let $\mathcal{C}$ be a specific set of $24$ pairwise disjoint $3$-cells of $P$ coming from an inscribed $24$-cell in $P$.  Let $\Gamma< \text{Isom}(\mathbb{H}^{4})$ be the discrete infinite covolume group generated by the reflections in the $96$ walls of $P$ that are not in $\mathcal{C}$. This paper shows that the inclusion map $\Gamma \hookrightarrow \text{Isom}(\mathbb{H}^{4})$ is infinitesimally rigid in the representation variety $\text{Hom}(\Gamma, \text{Isom}(\mathbb{H}^{4}))$, verifying a rigidity conjecture of Kerckhoff and Storm. 
Kerckhoff and Storm conjectured that compact hyperbolic $n$-orbifolds with totally geodesic boundary are infinitesimally rigid when $n>3$.  This paper verifies this conjecture for a specific example based on the $4$-dimensional hyperbolic $120$-cell.
\end{abstract}
\thanks{Aougab and Storm were partially supported by NSF grant DMS-0741604.}
\maketitle

\section{introduction}\label{introduction}

Given a discrete subgroup $\Gamma$ of a semisimple Lie group $G$, mathematicians have long studied the question of when $\Gamma$ can be deformed inside $G$.  On one side lie the great rigidity theorems of Calabi, Weil, Mostow, and Margulis, which roughly state that lattices in most semisimple Lie groups have no deformations.  On the other side lie the beautiful deformations theories of discrete groups of $2\times 2$ matrices, famous exceptions to the rigidity theorems.  If $G$ is a Lie group that is not represented by $2 \times 2$ matrices, and $\Gamma$ is not a lattice, then nearly nothing is known about the possible deformations of $\Gamma$.  Do reasonable geometric conditions exist that guarantee rigidity, or flexibility?

As a first step in this direction, Kerckhoff and Storm studied deformations of a specific discrete subgroup of $\text{Isom}(\mathbb{H}^4)$, the isometry group of hyperbolic $4$-space \cite{KS}.  This led them to the following rigidity conjecture.

\begin{conj}\label{KS rigidity}
If discrete group $\Gamma < \text{Isom}(\mathbb{H}^n)$ has Fuchsian ends, is convex cocompact, and $n>3$, then the inclusion map of $\Gamma$ is infinitesimally rigid.
\end{conj}

Qualitatively, groups $\Gamma$ satisfying the hypotheses of Conjecture \ref{KS rigidity} are very close to lattices.  Drop any of the three conditions on $\Gamma$, and the conjecture becomes false.  For any $n$, there are many convex cocompact discrete groups in $\text{Isom}(\mathbb{H}^n)$ that are not rigid, for example if $\Gamma$ is a free group.  For $n<4$ there exists a rich deformation theory applicable to discrete subgroups with Fuchsian ends.  Even weakening the convex cocompact condition to geometric finiteness makes the conjecture false.  A counterexample is presented in \cite{KS}.  Nonetheless, it seems reasonable to expect Conjecture \ref{KS rigidity} to be true.

Before trying to prove the full conjecture, it seems prudent first to verify it in a nontrivial case.  This is our goal here.  More specifically, we will construct an explicit reflection group in $\text{Isom}(\mathbb{H}^4)$ satisfying the hypotheses of Conjecture \ref{KS rigidity}, and verify that it is infinitesimally rigid.  The reflection group will be an infinite index subgroup of the lattice generated by reflections in the $3$-cells of the right-angled hyperbolic $120$-cell in $\mathbb{H}^4$.  This discrete group was chosen carefully to make the computations as simple as possible.  The group is described precisely in Sections \ref{120-cell} and \ref{TMT}.

The authors thank Daniel Allcock for his help with this research.

\section{Preliminaries}\label{sec:prelim}

This section is a brisk introduction to a few necessary facts from hyperbolic geometry.  For a detailed introduction see, for example, Ratcliffe's book \cite{Rat}.

Recall that Minkowski $(n+1)$-space, denoted as $\mathbb{M}^{n+1}$ is a real $(n+1)$-dimensional vector space equipped with the nondegenerate symmetric bilinear form 
\[ \left( \vec{x} , \vec{y} \right) := - x_0 y_0 + \sum_{i=1}^{n} x_i y_i.\]
A vector with positive norm is a space-like vector, a vector of norm $0$ is light-like, and a vector with negative norm time-like.  The hyperboloid model of hyperbolic $n$-space $\mathbb{H}^n$ is simply the the set of points
\[ \left\{ x \in \mathbb{M}^{n+1} \ | \  (x,x)= -1 \text{ and } x_0 > 0 \right\},\]
with the metric induced by $\mathbb{M}^{n+1}$.  Let $\text{O}(1,n)$ be the group of real $(n+1) \times (n+1)$ matrices $A$ such that $A^{T}MA= M$, where $M$ is the $(n+1) \times (n+1)$ diagonal matrix with diagonal entries $\{ -1,1,1,\ldots,1\}$.
\begin{comment}
\[\left( \begin{array}{rrrrr}
-1 & 0 & 0 & 0 & 0 \\      
0 & 1 & 0 & 0 & 0 \\
0 & 0 & 1 & 0 & 0 \\
0 & 0 & 0 & 1 & 0 \\
0 & 0 & 0 & 0 & 1
\end{array}\right).\]
\end{comment}
In other words, $\text{O}(1,n)$ is the group of linear transformations preserving the bilinear form $\left( \cdot, \cdot \right)$.
Then the group $\text{Isom}(\mathbb{H}^n)$ of isometries of $\mathbb{H}^n$ is the open subgroup of $\text{O}(1,n)$ given by matrices preserving the hyperboloid $\mathbb{H}^n$.  Note that orientation reversing isometries are included here; $\text{Isom}(\mathbb{H}^n)$ has two connected components.  Throughout the paper we will think of $\text{Isom}(\mathbb{H}^n)$ as this explicit group of matrices.  To simplify notation, define $G_n := \text{Isom}(\mathbb{H}^n)$.

Consider a discrete finitely generated subgroup $\Gamma$ of $G_n$ with presentation
\[ \Gamma = \left\langle r_1, r_2, \ldots, r_M \ | \ w_1 = w_2 = \ldots = w_N = 1 \right\rangle,\]
where each $w_i$ is a word in the generators $r_j$.  (The letter $r$ stands for ``reflection'', which will soon be the focus.)  To avoid degenerate cases, let us assume that $\Gamma$ is not virtually abelian, and does not preserve a copy of $\mathbb{H}^m$ for $m<n$.  If $\Gamma$ is torsion-free, then the quotient $\mathbb{H}^n/ \Gamma$ will be a hyperbolic manifold.  More generally, the quotient will be an orbifold.  If the quotient space has finite volume then $\Gamma$ is called a lattice.  Here we will be mainly interested in infinite volume quotient spaces, where a little more terminology is required.

\begin{defn}\label{fuchsian ends}
The discrete group $\Gamma$ has \emph{Fuchsian ends} if there exists a closed $\Gamma$-invariant convex set $C_\Gamma \subseteq \mathbb{H}^n$ that is an $n$-manifold with nonempty totally geodesic boundary, and such that the quotient $C_\Gamma / \Gamma$ has finite volume.
\end{defn}

Note that the boundary of $C_\Gamma$ will not be connected.  In the setting of groups with Fuchsian ends, $\Gamma$ is convex cocompact if and only if the quotient $C_\Gamma / \Gamma$ is compact. 

%If the quotient $C_\Gamma / \Gamma$ is compact then $\Gamma$ is said to be convex cocompact.

If $\Gamma$ has Fuchsian ends, then $\partial C_\Gamma$ consists of an infinite number of disjoint totally geodesic hyperplanes, each isometric to $\mathbb{H}^{n-1}$.  When $\Gamma$ is torsion-free, it has Fuchsian ends if and only if the quotient $C_\Gamma / \Gamma$ is a finite volume hyperbolic $n$-manifold with nonempty totally geodesic boundary.  In general, $C_\Gamma / \Gamma$ will be a finite volume orbifold with totally geodesic boundary.  For the reader familiar with convex cores, we note that if $\Gamma$ has Fuchsian ends then $C_\Gamma$ is the convex core of $\Gamma$.

Consider the representation variety $\text{Hom}(\Gamma,G_n)$.  The slightly larger space $\text{Hom}(\Gamma,\text{O}(1,n))$ can be represented explicitly as the zero set of a collection of real polynomials as follows.  A homomorphism $\rho: \Gamma \rightarrow \text{O}(1,n)$ is determined uniquely by the $M$ matrices $\rho(r_j)$ of $\text{O}(1,n)$, which can be thought of as a point in a real vector space $V$ of dimension $M \cdot (n+1)^2$.  The $\rho(r_j)$ must lie in $\text{O}(1,n)$, meaning for each $j$ the $(n+1)^2$ polynomials $\rho(r_j)^{T} M \rho(r_j) = M$ are satisfied.   Each relation $w_i$ of $\Gamma$ becomes a system of $(n+1)^2$ polynomial equations in the entries of the matrices $\rho(r_j)$.  The variety $\text{Hom}(\Gamma,\text{O}(1,n)) \subset V$ is precisely the set of points where these $(M+N) \cdot (n+1)^2$ polynomial equations $\{ P_\alpha\}$
are satisfied.  Finally, we are interested in the subset $\text{Hom}(\Gamma,G_n) \subset \text{Hom}(\Gamma,\text{O}(1,n))$.  It is the set of connected components where the upper-left matrix entry of each $\rho(r_j)$ is positive.  (This entry cannot be zero for a matrix in $\text{O}(1,4)$.)

The group $G_n$ acts on $\text{Hom}(\Gamma,G_n)$ by conjugation as follows:
\[ (g \cdot \rho) (\gamma) := g \rho(\gamma) g^{-1},\]
where $g\in G_n$ and $\gamma \in \Gamma$.  It is easy to see that this action is algebraic.
The inclusion map $\Gamma \rightarrow G_n$ is a point $\iota \in \text{Hom}(\Gamma,G_n)$.  In general, analyzing the orbits of this $G_n$-action can be complicated.  However, using our assumptions that $\Gamma$ is not virtually abelian and does not preserve a lower dimensional hyperbolic space, it is easy to show that the orbit $G_n \cdot \iota \subseteq \text{Hom}(\Gamma,G_n)$ is a manifold of dimension equal to that of $G_n$ \cite{Gold1}.

A homomorphism $\rho \in \text{Hom}(\Gamma,G_n)$ is \emph{locally rigid} if it has an open neighborhood contained inside the orbit $G_n \cdot \rho$, in other words all nearby representations are obtained by conjugation.  The infinitesimal analogue of this notion is useful.  An \emph{infinitesimal deformation} of $\rho$ is a tangent vector $p'(0) \in T_\rho V$ to a smooth path $p : (-\varepsilon,\varepsilon) \rightarrow V$ such that 
\[ \frac{d}{dt} |_0 P_\alpha (p(t) ) = 0\]
for the $(M+N) \cdot (n+1)^2$ polynomial equations defining $\text{Hom}(\Gamma,\text{O}(1,n))$ (and locally defining $\text{Hom}(\Gamma,G_n)$).
This slightly odd definition is necessitated by the possibility that $\text{Hom}(\Gamma,G_n)$ is singular.  Similarly, an \emph{infinitesimal conjugation} of $\rho$ is a tangent vector $p'(0) \in T_\rho V$ to a smooth path $p$ in the orbit $G_n \cdot \rho$.  We say $\rho$ is \emph{infinitesimally rigid} if every infinitesimal deformation is an infinitesimal conjugation.  An infinitesimally rigid homomorphism is also locally rigid \cite{We}.  By the above discussion, to show that the inclusion map $\iota$ of $G_n$ is infinitesimally rigid, it suffices to show that the linear subspace of infinitesimal deformations has dimension equal to that of $G_n$.

When $\Gamma < G_n$ is a lattice, we have the famous rigidity theorem of Calabi, Weil, and Garland.

\begin{thm}\label{Weil rigidity}\cite{Cal,We,Garland}
If $\Gamma < \text{Isom}(\mathbb{H}^n)$ is a lattice and $n>3$ then the inclusion map of $\Gamma$ is infinitesimally rigid.
\end{thm}

For discrete groups that are not lattices, Theorem \ref{Weil rigidity} is false, but it is interesting to study whether or not a similar rigidity theorem might hold for any other natural class of discrete groups.  Looking for infinite covolume groups which are as ``close'' to lattices as possible, Kerckhoff and Storm were led to consider discrete groups with Fuchsian ends.  In $\mathbb{H}^3$, groups with Fuchsian ends are not rigid, and have a beautiful deformation theory investigated by many people \cite{Th}.  In higher dimensions, simple naive dimension counts suggest that groups with Fuchsian ends are rigid.  These dimension counts are, of course, not rigorous.  Nonetheless, Kerckhoff and Storm were led to the following conjecture, stated first in the introduction.

\begin{named}{Conjecture \ref{KS rigidity}}
If discrete group $\Gamma < \text{Isom}(\mathbb{H}^n)$ has Fuchsian ends, is convex cocompact, and $n>3$, then the inclusion map of $\Gamma$ is infinitesimally rigid.
\end{named}

Note that Conjecture \ref{KS rigidity} is false without the assumption that $\Gamma$ is convex cocompact.  A counterexample was studied in \cite{KS}.  As a first step toward proving this conjecture, the goal of this paper is to verify it in a specific $4$-dimensional example.  In search of an example where the computations are as simple as possible, the first place to look is among hyperbolic reflection groups, which we now introduce briefly.  For more information the authors recommend \cite{Vin2}.

A reflection isometry in $G_4 = \text{Isom}(\mathbb{H}^4)$ is given by a matrix of $\text{O}(1,4)$ with a $4$-dimensional eigenspace of eigenvalue $1$, and a single eigenvalue equal to $-1$ whose corresponding eigenvector is space-like.  It fixes a codimension $1$ totally geodesic hyperplane of $\mathbb{H}^4$ given by the intersection of its $+1$-eigenspace with $\mathbb{H}^4$.  The orthogonal complement of the fixed hyperplane is reflected across the hyperplane by the isometry.  The reflection isometry is determined uniquely by the hyperplane and vice versa.

Consider the interesting special case where $\Gamma < G_4$ is a group generated by $M$ reflections $r_1$, $r_2$, $\ldots$, $r_M$ with corresponding hyperplanes $H_1$, $H_2$, $\ldots$, $H_M$, and the hyperplanes bound a convex (possibly infinite volume) region $P$ of $\mathbb{H}^4$ known as a fundamental domain.  Moreover, assume the dihedral angles of $P$ are all $\pi/2$.  In this case the presentation of $\Gamma$ takes the following nice form \cite{Vin2}:
\[\Gamma = \left\langle r_1, r_2, ...,r_M \ | \ (r_ir_j)^{m_{ij}}= r_i^2 = 1  \right\rangle ,\]
where for each pair $1 \le i < j \le M$, $m_{ij}$ is either $2$ or infinity.  As usual, the ``relation'' $(r_i r_j)^\infty = 1$ indicates the product $r_i r_j$ has infinite order.  The product of reflections $r_i r_j$ has order $2$ when the hyperplanes $H_i$ and $H_j$ intersect at angle $\pi / 2$.  This is best seen by thinking about the picture in the plane of the product of two reflections.  Otherwise the product has infinite order, or $m_{ij} = \infty$.  (Here we are not allowing intersections at other angles.)

We will be considering the representation variety $\text{Hom}(\Gamma,G_4) \subset \mathbb{R}^{25 M}$ for groups $\Gamma$ of the above form.

The following lemma shows that near the inclusion map in $\text{Hom}(\Gamma,G_4)$ any representation has the property that it maps the generators $r_i$ to reflection isometries of $G_4$.

\begin{lem} \label{reflection lem}
Let $A \in G_n$ be a reflection isometry. Then $A$ has an open neighborhood $U\subset G_n$ such that if $B \in U$ and $B^{2}= \text{Id}$, then $B$ is a reflection isometry. 
\end{lem}

\begin{pf}
Recall that any involution can be diagonalized over $\mathbb{C}$ such that the diagonal entries are all $\pm1$.  Recall also that a reflection must have exactly one eigenvalue equal to $-1$.  We can define the open neighborhood $U$ of $A$ to be all matrices in $G_n$ having $(n-1)$ positive eigenvalues (not necessarily distinct) and $1$ negative eigenvalue. Then any matrix in $U$ that is an involution will necessarily be a nontrivial isometry which fixes a hyperplane of codimension $1$, in other words, a reflection.
\end{pf}

\begin{cor} \label{cor reflection lem}
 There exists an open neighborhood $U \subset \text{Hom}(\Gamma, G_4)\subset \mathbb{R}^{25 M}$ of the inclusion map such that if $\rho \in U$, then $\rho(r_{i})$ is a reflection for all $i$. 
 \end{cor}

The representation variety $\text{Hom}(\Gamma, G_4)$ sits naturally as (connected components of) a subvariety of $\mathbb{R}^{25 M}$, but by exploiting the fact that locally all the generators are reflections, it is possible to reduce the dimension considerably.  To a reflection isometry $r$ with fixed hyperplane $H$ we can associate the $1$-dimensional subspace of $\mathbb{M}^5$ given by vectors Minkowski orthogonal to $H$.  This subspace is simply the $+1$-eigenspace of the matrix corresponding to $r$, and it will consist of space-like vectors.  In reverse, choosing a space-like vector $\vec{n}$ determines a codimension $1$ hyperplane $\vec{n}^\perp \cap \mathbb{H}^4$, which in turn determines a reflection isometry.  Letting $\mathcal{S} \subset \mathbb{M}^5$ denote the set of space-like vectors, this process defines a map
\[ \mu: \mathcal{S}^M \rightarrow G_4^M \subset \mathbb{R}^{25 M} \]
with image equal to the set of $M$-tuples of reflection isometries.  

Moreover, the geometry of the fundamental domain $P$ can be read from the normal vectors.  Specifically,
let $\vec{n}_{1}$ and $\vec{n}_{2}$ be space-like in $\mathbb{M}^5$.  Then let 
\[ H_{i}= \left\{ \vec{v}\in \mathbb{M}^5 \ | \  \left(\vec{v},\vec{n}_{i}\right) = 0 \right\}.\]  
If $H_{1}\cap H_{2}\cap\ \mathbb{H}^{4} \neq \emptyset$ and $H_{1}$ and $H_{2}$ intersect at angle $\theta$, then 
\[ \frac{\left(\vec{n}_{1},\vec{n}_{2}\right)}{\sqrt{\left(\vec{n}_{1}, \vec{n}_{1}\right)\left(\vec{n}_{2}, \vec{n}_{2}\right)}} = -\cos\left(\theta\right). \] 
In particular, we will consider only examples where $\theta$ is always $\pi/2$, in which case the above simplifies to
\[\left(\vec{n}_{1},\vec{n}_{2}\right) = 0. \]

Choose a representation $\rho \in \text{Hom}(\Gamma,G_4)$ sending each generator $r_i$ to a reflection.  In the above manner, we replace every reflection matrix $\rho(r_i)$ by a corresponding normal vector $\vec{n}_i$.  
We now have $M$ Minkowski $5$-vectors, instead of $M$ $5 \times 5$ matrices.  The relations of $\Gamma$ written in terms of normal vectors all have the form:
 \[\left(\vec{n}_{i}, \vec{n}_{j}\right)= 0 \text{ when } m_{ij} =2 \text{ (angle condition)}\]

Let $U \subset \text{Hom}(\Gamma,G_4)$ be the open neighborhood of the inclusion map of Corollary \ref{cor reflection lem}.  It is clear that, restricted to $\mu^{-1}(U)$, $\mu$ is a submersion with fibers given by scaling the normal vectors $\vec{n}_i$.  Given a vector $v = \left( \dot{\vec{n}}_i \right)$ tangent to the point $( \vec{n}_i ) \in \mu^{-1}(\rho)$, $\mu_* v \in T_\rho \mathbb{R}^{25 M}$ is an infinitesimal deformation of $\rho$ if and only if $v$ satisfies the following derivatives of the above angle conditions:  
\[ \left(\dot{\vec{n}}_{i}, \vec{n}_{j}\right)+ \left(\vec{n}_{i}, \dot{\vec{n}}_{j}\right)= 0 \text{ when } m_{ij} =2\]
Any solution to this system of polynomials will now be a vector $( \dot{\vec{n}}_i )$ in $\mathbb{R}^{5 M}$ instead of $\mathbb{R}^{25 M}$. 
Note that we do not need to include the relation corresponding to each generator being an involution.  This interpretation clearly has an immense computational advantage over working with the matrices directly.

Finally, to show $\rho$ is infinitesimally rigid in $\text{Hom}(\Gamma,G_4)$ it suffices to show that the subspace of vectors $( \dot{\vec{n}}_i ) \subset \mathbb{R}^{5 M}$ satisfying the above equations has dimension equal to $10 + M$, where $10 = \text{dim}(G_4)$ dimensions come from infinitesimal conjugations by $G_4$, and there are $M$ dimensions corresponding to scaling the normal vectors, one for each generator of $\Gamma$. 

\section{The $120$-cell}\label{120-cell}
Here we will outline some of the basic properties of the hyperbolic $120$-cell.  We will use the word face to indicate a $2$-cell, and wall to indicate a $3$-cell. 
The $120$-cell is a regular convex polytope formed by $120$ dodecahedral walls, where we define regular to mean that its symmetry group acts transitively on the set of flags.  (A flag of the $120$-cell is quadruple consisting of $1$ point, $1$ edge containing the point, $1$ pentagonal face containing the edge, and $1$ dodecahedral wall containing the face.)  The $120$-cell has $600$ vertices, and is dual to the $600$-cell formed by $600$ icosahedra.

The $120$-cell can be embedded into $\mathbb{H}^4$ in such a way that the resulting convex hyperbolic polytope is compact, regular, and all intersecting dodecahedra hit at right angles.  The quickest way to describe this polytope is by giving a list of $120$ space-like vectors that are Minkowski-normal to the $120$ walls.  For completeness, we list the $120$ normal vectors in $\mathbb{M}^5$, using the Golden ratio $\tau = (1+\sqrt{5})/2$: \cite{Cox4}
\begin{enumerate}
\item\label{1}  The $8$ vectors obtained by permuting the last $4$ coordinates of $\left(\sqrt{2\tau}, \pm 2, 0, 0, 0\right)$.
\item\label{2}  The $16$ vectors of the form $\left(\sqrt{2\tau},\pm1,\pm1,\pm1,\pm1\right)$.
\item  The $96$ even permutations in the last $4$ coordinates of $\left(\sqrt{2 \tau}, \pm\tau, \pm1, \pm\tau^{-1}, 0\right)$.
\end{enumerate}
Consider the set $\mathcal{C}$ of walls of the $120$-cell given by the $24$ space-like vectors of items (\ref{1}) and (\ref{2}) above.  Interestingly, these $24$ vectors determine $24$ hyperplanes of $\mathbb{H}^4$ which bound a regular (hyperideal) hyperbolic $24$-cell, which is a polytope with $24$ octahedral walls.  This set $\mathcal{C}$ will play an important role here.  
One can compute directly that the walls of $\mathcal{C}$ are pairwise disjoint in $\mathbb{H}^4$.
Moreover, there does not exist a set of $25$ pairwise disjoint walls of the $120$ cell.

\begin{prop} \label{24 walls}
Let $P$ be the hyperbolic $120$-cell.  A maximum set of pairwise disjoint walls of $P$ has 24 elements. 
\end{prop}

\begin{pf}
Suppose we remove $24$ walls from the $120$-cell. For each wall that we remove, we place a marker on any wall which was adjacent to the removed wall. Each wall of the $120$-cell is adjacent to $12$ other walls; therefore after removing $24$ walls we have placed $24\times12= 288$ markers on the remaining $96$ walls, or an average of $3$ markers per wall that remains. The claim is that every remaining wall has exactly $3$ markers. If this is indeed true, then removing an additional $25^{\text{th}}$ wall would be impossible because any wall we removed would be adjacent to $3$ of the walls in the original set of $24$.

Suppose one of the walls had only $2$ markers.  This implies that there is at least one other wall with $4$ markers, implying that the intersection of a certain set of $4$ pairwise disjoint walls with this wall of $4$ markers is a set of $4$ pairwise disjoint pentagonal faces of a dodecahedron. This is impossible because the maximum number of pairwise disjoint faces of a dodecahedron is $3$. Therefore no wall has $2$ or fewer markers implying every wall has $3$.
\end{pf}
\noindent The above proof was explained to us by Allcock.

\section{The Main Theorem} \label{TMT}

We are now ready to discuss the main result. 

\begin{thm} \label{Main Thm}
Let $P \subset\ \mathbb{H}^{4}$ be the hyperbolic regular right-angled $120$-cell. Let $\mathcal{C}$ be the set of $24$ pairwise disjoint walls of $P$ given in Section \ref{120-cell}.  Let $\Gamma < G_4 = \text{Isom}(\mathbb{H}^{4}) < \text{O}(1,4)$ be the discrete infinite covolume group generated by the reflections in the $96$ walls of $P$ that are not in $\mathcal{C}$.  Then the inclusion map $\Gamma \rightarrow G_4$ is infinitesimally rigid in the representation variety $Hom(\Gamma, G_4)$.
\end{thm} 

Before beginning the proof let us discuss the connection to Conjecture \ref{KS rigidity}.  

\begin{prop} \label{Fuchsian end prop}
Let $\Gamma < G_4$ be the discrete group of Theorem \ref{Main Thm}.  Then $\Gamma$ has Fuchsian ends and the quotient $C_\Gamma / \Gamma$ is isometric to the $120$-cell $P$.  In particular, $\Gamma$ is convex cocompact.
\end{prop}

\begin{pf}
Consider the collection of $24$ dodecahedral walls $\mathcal{C} \subset \mathbb{H}^4$ and its orbit $\Gamma \cdot \mathcal{C}$.  Using the facts that all the dihedral angles of $P$ are $\pi/2$, and all the walls of $\mathcal{C}$ are pairwise disjoint, it follows that any intersecting translates of $\mathcal{C}$ in the orbit $\Gamma \cdot \mathcal{C}$ glue together smoothly.  In particular, the entire orbit is a disjoint union of totally geodesic hyperplanes.  Let $C_\Gamma \subset \mathbb{H}^4$ be the $\Gamma$-invariant convex subset bounded by $\Gamma \cdot \mathcal{C}$.  This shows that $\Gamma$ has Fuchsian ends.  It is clear that the quotient $C_\Gamma / \Gamma$ is isometric to $P$.
\end{pf}

Combining Theorem \ref{Main Thm} and Proposition \ref{Fuchsian end prop} shows that Conjecture \ref{KS rigidity} is true for $\Gamma$, which we record here.

\begin{cor} 
The discrete subgroup $\Gamma < G_4$ has Fuchsian ends, is convex cocompact, and infinitesimally rigid, as predicted by Conjecture \ref{KS rigidity}.
\end{cor}

Let us begin the proof of Theorem \ref{Main Thm}.  We begin with $96$ $5 \times 5$ reflection matrices in $G_4<\text{O}(1,4)$, the image of the generators of $\Gamma$ under the inclusion map.   We can replace each such reflection matrix with one of its (space-like) eigenvectors $\vec{n}_i$ corresponding to eigenvalue $-1$.  (All the possible choices are colinear.)  The complete list of $96$ vectors $\vec{n}_i$ will be 
the $96$ even permutations in the last $4$ coordinates of $\left(\sqrt{2 \tau}, \pm\tau, \pm1, \pm\tau^{-1}, 0\right)$, where $\tau$ is the Golden ratio.  Pairs $\vec{n}_i$ and $\vec{n}_j$ corresponding to orthogonal walls of $P$ will satsify the relation
\[ \left( \vec{n}_i , \vec{n}_j \right) = 0.\]

Let $W \subset T_{\vec{n}} \mathbb{R}^{480}$ be the linear space of infinitesimal deformations.  Suppose we are given an infinitesimal deformation 
\[\left( \dot{\vec{n}}_i \right) \in W \subset T_{\vec{n}} \mathbb{R}^{480} \cong \mathbb{R}^{480}\]
of the $96$ space-like vectors which maps via $\mu_*$ to an infinitesimal deformation of $\Gamma$ in $G_4$.  Then the following linear equation must hold for all pairs $i,j$ corresponding to intersecting orthogonal walls of $P$:
\[ \left(\dot{\vec{n}}_{i}, \vec{n}_{j}\right)+ \left(\vec{n}_{i}, \dot{\vec{n}}_{j}\right)= 0.\]
Keep in mind that the $24$ walls of $\mathcal{C}$ have been removed from $P$.  There are $432$ equations of this type, corresponding to the $432$ remaining faces of $P$.
As described above, we wish to find the solution space $W$ of this set of linear polynomials.  To prove that $\Gamma$ is infinitesimally rigid it suffices to show that the solution space has $10 + 96$ dimensions, where $10$ comes from the action of $G_4$ by conjugation, and $96$ comes from scaling each of the $96$ space-like vectors $\vec{n}_i$.  These $96$ scaling dimensions are killed by $\mu_*$.  Indeed, by this count we see that the solution space $W$ has at least $106$ dimensions, so it remains to show that it has at most $106$ dimensions.

Begin by enumerating $\{ p_1, p_2, \ldots, p_{432} \}$ the polynomials defining $W$.  Each $p_i$ is a polynomial in the variables $(v_1, v_2, \ldots, v_{480})$ where 
\[\dot{\vec{n}}_{i}= \left(v_{5(i-1)+1}, v_{5(i-1)+2}, v_{5(i-1)+3}, v_{5(i-1)+4}, v_{5i}\right).\]
Then define the matrix $A_{\text{alg}}$ of algebraic numbers
\[ (A_{\text{alg}})_{ij} := \frac{d}{dv_j} \left( p_i (v_1, v_2, \ldots, v_{480} ) \right).\]
Each entry $(A_{\text{alg}})_{ij}$ is a number because each polynomial $p_i$ is linear in the variables $v_j$.  $A_{\text{alg}}$ has $432$ rows and $480$ columns.  The goal is to show its kernel has dimension at most $106$.  We will do this by showing the rank is at least $374$.

In order to ensure that Mathematica can correctly and efficiently run the computation, we replace every algebraic number entry of $A_{\text{alg}}$ with a square matrix of rational numbers. The number field $\mathbb{Q}(\alpha)$ contains every entry of $A_{\text{alg}}$ when $\alpha = \sqrt{1 + \sqrt{5}}$.
We choose the basis 
\[ \left\{ 1, \alpha, \alpha^2, \alpha^3 \right\}\]
for $\mathbb{Q}(\alpha)$ as a vector space over $\mathbb{Q}$.
Then each entry of $A_{\text{alg}}$ is replaced by a $4 \times 4$ rational matrix representing its action by left multiplication on $\mathbb{Q}(\alpha)$ with respect to this basis.  
For example, the number $\frac{1}{2}\left(-1-\sqrt{5}\right)$ is replaced by the $4\times4$ matrix
\[\left(\begin{array}{clrr} % 
     0 & 0 & -2 & 0 \\     
     0 & 0 & 0            & -2\\
     -\frac{1}{2} & 0 & -1           &0\\
      0& -\frac{1}{2} & 0            &-1    
       \end{array}\right)
      \] 

Once we have the rational $4\times4$ matrix replacements for each entry of $A_{\text{alg}}$, it becomes a $432 \times 480$ matrix of $4 \times 4$ matrices.  By ignoring the structure of the $4 \times 4$ matrices, we can think of it as a $1728\times1920$ matrix $A_{\mathbb{Q}}$ of rational numbers. 
It is easy to see that the kernel of $A_{\mathbb{Q}}$ has dimension at least $4$ times that of $A_{\text{alg}}$.  (If $u$ is in the kernel of $A_{\text{alg}}$ then $\{ u,\alpha u, \alpha^2 u, \alpha^3 u \}$ will be in the kernel of $A_\mathbb{Q}$ and linearly independent over $\mathbb{Q}$.)  Therefore the rank of $A_\mathbb{Q}$ is at most $4$ times the rank of $A_{\text{alg}}$.  To prove the infinitesimal rigidity of $\Gamma$, it therefore suffices to show that the rank of $A_\mathbb{Q}$ is at least $4 \cdot 374 = 1496$.

In order to further simplify the computation for Mathematica, multiply each row of $A_\mathbb{Q}$ by the least common multiple of the denominators of the rational numbers in that row.  This gives us the integer matrix $A_{\mathbb{Z}}$ with the same rank as $A_\mathbb{Q}$.  With some effort on a powerful computer (in 2008), Mathematica was able to determine that the rank of this matrix is $1496$, the desired result.
Therefore $A_\text{alg}$ has a kernel of dimension at most 106.
This proves that the inclusion map of the discrete group $\Gamma < G_4$ is infinitesimally rigid in the representation variety $\text{Hom}(\Gamma,G_4)$.  

The rank calculation for $A_{\mathbb{Z}}$ took approximately $8$ hours, and several gigabytes of RAM.  This is perhaps related to the well known computational difficulty that diagonalizing an integer matrix via the Euclidean algorithm can lead to intermediate matrix entries of enormous size.  We can greatly decrease the calculation time by reducing $A_{\mathbb{Z}}$ modulo a suitable prime $p$, obtaining matrix $A_p$. When working in the finite field $\mathbb{F}_{p}$ the diagonalization process is nearly trivial, and $\text{rank}(A_p) \leq \text{rank}(A_\mathbb{Z})$.
Therefore, we simply reduce $A_{\mathbb{Z}}$ modulo several primes until we find one such that $\text{rank}(A_p)= 1496$.  Once we find such a prime we are done by the above inequality.  
Using the prime $113$, this method allows us to arrive at the correct answer of a rank of $1496$ for $A_{\mathbb{Z}}$ with a calculation time of approximately $5$ seconds.  

\subsubsection*{Acknowledgements:}

All computer calculations were performed in $2008$ using Wolfram's Mathematica $6.0.1.0$ on a MacPro running $64$-bit Linux.

\begin{comment}
\section{References}
1. Artin, Michael. Algebra. \vskip 6 pt
2. Coxeter, HSM. Regular Complex Polytopes. \vskip 6 pt
3. Coxeter, HSM. Regular Polytopes. \vskip 6 pt
4. Coxeter, HSM. Introduction to Geometry. \vskip 6 pt
5. Davis, Michael W. A Hyperbolic 4 Manifold \vskip 6 pt
6. Kerckhoff, Steven P. and Storm, Peter A. From the 24-cell to the Cuboctahedron \vskip 6 pt
7. Gallier, John. www.cis.upenn.edu/~cis610/cis610sl7.pdf \vskip 6 pt
8. Ratcliffe, John G. Foundations of Hyperbolic Manifolds.
\end{comment}

\bibliography{bibliography}
\end{document}